\documentclass[a4paper,12pt]{article}
\begin{document}

\title{Schauder bases in Dirac modules over quaternions.}
\author{Sergey V. Ludkovsky}
\date{19 October 2013}
\maketitle

\begin{abstract}
Dirac modules over the quaternion skew field are investigated on a
compact domain relative to the supremum norm and Hardy's norm with
the parameter $1<p<\infty $ as well. An existence of Schauder bases
in them is proved. Procedures for construction of such bases are
outlined.

\end{abstract}

\section{Introduction.}
The theory of Schauder bases in Banach spaces is the important part of Functional Analysis \cite{gurlusb,jarchb,lindliorb,woytb}.
As it is known now it exists not in all Banach spaces.
But an existence or finding of it in concrete classes of Banach spaces is frequently a serious problem
related with their particular structure \cite{gurlusb,jarchb,lanc92,lanc97,lindliorb,lusky92}-\cite{lusky04,woytb}.
\par On the other hand, hypercomplex analysis over Clifford algebras
and quaternions in particular is developing fast (see \cite{brdeso},
\cite{gentstram07} - \cite{guespr},
\cite{ludoyst,ludfov,lufscdvm,lufjmsrf}, \cite{sudbcps79} and
references therein). It has many applications not only in
mathematics, but also in natural sciences
\cite{gilmurr,girard,guetze}. For example it permits to integrate
new types of partial differential equations (see
\cite{guespr,ludancdb} and references therein).
\par One of the classical examples is the solution of Klein-Gordon's hyperbolic partial differential equation
by Dirac with the help of quaternions \cite{dirac,gilmurr,girard}.
Since that time the partial differential operator of the first order
over quaternions used by Dirac for the decomposition of the
hyperbolic partial differential Klein-Gordon operator is widely used
and is frequently known under his name.
\par The quaternion skew field $\bf H$ is associative and non-commutative.
It is the algebra over the real field $\bf R$, but it is not the algebra
over the complex field $\bf C$, because the center of $\bf H$ is the real field.
\par Each complex holomorphic or harmonic function has locally a
power series expansion, but for their concrete Banach spaces it is
frequently a serious problem whether they have Schauder bases
\cite{bock74,lanc92,lanc97,lusky92}-\cite{lusky03,woytb}.
\par In the paper \cite{bock74} an existence of a Schauder basis in the Banach algebra $A(K)$ of all
holomorphic functions on the open unit circle $K$ in the complex
field $\bf C$ having continuous bounded extensions on its closure
$\bar K$ was proved and a procedure for its construction was
described. \par Characterizing complex holomorphic functions by the
condition $\bar {\partial }f(x)=0$ one can consider their quaternion
analog $\sigma f=0$, where $\bar {\partial }f(x)=\partial
f(x)/\partial {\bar x}$ for the complex variable $x$, whilst $\sigma
$ is the Dirac operator over the quaternion skew field.
Nevertheless, the technique presented in \cite{bock74} does not work
over the quaternion skew field $\bf H$ because of specific features
of $\bf H$ and $ker (\sigma )$. Moreover, in that work very
particular properties of complex holomorphic functions on $K$,
functions on the circumference $S^1$ and the commutativity of the
complex field were used for the proof, which are not valid over the
quaternion skew field. For example, converging power series over
quaternions are more complicated and contain additives like
$a_{n,1}z^{k_1}...a_{n,m}z^{k_m}$, where $z$ is a variable in a
domain $\Omega $ in $\bf H$, $k_1,...,k_m$ are non-negative
integers, $m$ is a natural number, while $a_{n,1},...,a_{n,m}$ are
quaternion coefficients. Moreover, in \cite{bock74} it was used that
each complex number $x$ can be written as $x=|x|e^{i\phi }$, where
$\phi $ is a real parameter (argument of $x$). This is not the case
in $\bf H$, because $z=|z| e^M$ with $M$ in the non-commutative
purely imaginary domain ${\cal I} = {\bf R}i_1\oplus {\bf
R}i_2\oplus {\bf R}i_3\subset \bf H$ such that $Re (M)=0$, but
generally $e^M$ and $e^N$ do not commute, when purely imaginary
quaternions $M, N\in {\cal I}$ do not commute. \par Using the
commuting bounded approximation property it was proved anew in
\cite{lusky92} that the disk algebra $A(K)$ has a basis and that the
Hardy space $H^p(K)$ has an unconditional basis. On the
multi-dimensional torus $K^N$ these results were extended in
\cite{lanc92,lanc97}.
\par In this work the kernel of the Dirac operator $\sigma $ is investigated on a domain $\bar {\Omega }$ in $\bf H$ quasi-conformal with the closed unit ball $\bar B$ (see Section 3.1). The Dirac operator is considered from the space of continuously differentiable functions $f$ on $\bar {\Omega }$ into the space of continuous functions on $\bar {\Omega }$. That is $\sigma f$ is taken on an open domain $\Omega $ such that the function
$f$ has a continuous bounded extension on $\bar{\Omega }$, where
$\bar{\Omega }$ is the closure of $\Omega $. It appears to be a left
module over the quaternion skew field. But apart from the complex
case it is not an algebra even over the real field. In our paper we
elaborate a new technique different from previous works.
\par In Section 3 of this article Dirac modules over the quaternion skew field $\bf H$ are
investigated. For this purpose an analog of the Stone-Weierstrass
theorem over quaternions is proved in Section 2. An existence of a
Schauder basis in a Dirac module supplied with the supremum norm is
proved in Theorem 3.15 and Corollary 3.16. Procedures for
construction of such bases are outlined. Moreover, relative to
Hardy's norm with the parameter $1<p<\infty $ Theorem 3.18 about
unconditional bases in Dirac modules is proved.
\par Main results of this paper are obtained for the first time.

\section{Banach spaces over the quaternion skew field.}
{\bf 1. Definitions and Notes.} An ${\bf R}$ linear space $X$ which
is also left and right ${\bf H}$ module will be called an ${\bf H}$
vector space. We present $X$ as the direct sum
\par $(DS)$ $\quad X=X_0i_0\oplus ... \oplus X_3 i_3$, where
$X_0$,...,$X_3$ are pairwise isomorphic real linear spaces, where
$i_0,...,i_3$ are generators of the quaternion skew field ${\bf H}$
such that $i_0=1$, $i_k^2=-1$ and $i_ki_j=-i_ji_k$ for each $k\ge 1$
and $j\ge 1$ so that $k\ne j$.
\par  Let $X$ be an
$\bf R$ linear normed space which is also left and right ${\bf H}$
module such that  \par $(1)$ $0\le \| ax \|_X = |a| \| x \|_X $ and
\par $(2)$ $ \| x a \|_X= |a| \| x \|_X $ and \par $(3)$ $\| x+y \|_X \le \| x \|_X + \| y \|_X $ \\
for all $x, y\in X$ and $a\in {\bf H}$. Such space $X$ will be
called an ${\bf H}$ normed space.
\par Suppose that $X$ and $Y$ are two normed spaces over ${\bf H}$.
A continuous $\bf R$ linear mapping $\theta : X\to Y$ is called an
$\bf R$ linear homomorphism. If in addition $\theta (bx)=b\theta
(x)$ and $\theta (xb)=\theta (x)b$ for each $b \in {\bf H}$ and
$x\in X$, then $\theta $ is called a homomorphism of ${\bf H}$ (two
sided) modules $X$ and $Y$.
\par If a homomorphism is injective, then it is called an embedding
($\bf R$ linear or for ${\bf H}$ modules correspondingly).
\par If a homomorphism $h$ is bijective and from $X$ onto $Y$ so that its inverse
mapping $h^{-1}$ is also continuous, then it is called an
isomorphism ($\bf R$ linear or of ${\bf H}$ modules respectively).

\par {\bf 2. Definitions.} One says that a real vector space $Z$ is supplied with a scalar
product if a bi-$\bf R$-linear bi-additive mapping $<,>: Z^2\to {\bf
R}$ is given satisfying the conditions:
\par $(1)$ $<x,x> ~ \ge 0$, $ ~ <x,x>=0$ if and only if $x=0$;
\par $(2)$ $<x,y>=<y,x>$;
\par $(3)$ $<ax+by,z>=a<x,z>+b<y,z>$ for each real numbers $a, b\in {\bf R}$ and vectors $x, y,
z\in Z$.
\par Then an ${\bf H}$ vector space $X$ is supplied with an
${\bf H}$ valued scalar product, if a bi-${\bf R}$-linear bi-${\bf
H}$-additive mapping $<*,*>: X^2\to {\bf H}$ is given such that
\par $(4)$ $\quad <f,g> = \sum_{j,k} <f_j,g_k>i_j^*i_k$, \\
where $f=f_0i_0+...+f_3i_3$, $ ~f, g\in X$, $ ~ f_j, g_j \in X_j$,
each $X_j$ is a real linear space with a real valued scalar product,
$(X_j, <*,*>)$ is real linear isomorphic with $(X_k, <*,*>)$ and
$<f_j,g_k>\in {\bf R}$ for each $j, k$. The scalar product induces
the norm: \par $(5)$ $\| f \| := \sqrt{<f,f>}$.
\par An ${\bf H}$ normed space or an
${\bf H}$ vector space with an ${\bf H}$ valued scalar product
complete relative to its norm will be called an ${\bf H}$ Banach
space or an ${\bf H}$ Hilbert space respectively.

\par {\bf 3. Banach spaces of continuous functions over the quaternion skew field.}
As usually the quaternion skew field ${\bf H}$ is supplied with its
standard norm topology: $|z|=\sqrt{zz^*}$ for each $z\in {\bf H}$.
Considered as the real normed space the quaternion skew field ${\bf
H}$ has the real shadow which is the Euclidean space ${\bf R}^4$.
\par Let $C(U,{\bf H})$ denote the set of all continuous
functions on a canonical closed subset $U$ in ${\bf H}$.  It is an
$\bf R$-linear space and a left- and right- ${\bf H}$-module.
Moreover, $C(U,{\bf H})$ is the algebra over the real field with the
point-wise addition and multiplication of functions. Supply it with
the norm
\par $(1)$ $\| f \| _{C(U,{\bf H})} := \sup_{z\in U} |f(z)|$.
\par Relative to this norm the space $C(U,{\bf H})$ is Banach.
There are the following quaternion analogs of the Stone-Weierstrass
theorem.
\par {\bf 4. Theorem.} {\it Suppose that $A$ is a separating points
subalgebra in $C(T,{\bf H})$ such that $Z(A) = \emptyset $ (that is,
functions from $A$ have no any common zero in $T$), where $T$ is a
Hausdorff compact topological space, then this set $A$ is
(everywhere) dense in $C(T,{\bf H})$.}
\par {\bf Proof.} For proving this theorem we use the classical
Stone-Weierstrass theorem. It states: let $A$ be a separating points
subalgebra in $C(T,{\bf R})$ such that $Z(A)=\emptyset $, where $T$
is a Hausdorff compact topological space, then this set $A$ is
(everywhere) dense in $C(T,{\bf R})$ (see, for example, \S 4.10
Theorem A in \cite{edwardsb}). The algebra $C(U,{\bf H})$ is
isomorphic with $[C(U,{\bf R})i_0]\oplus ... \oplus [C(U,{\bf
R})i_3]$. This decomposition induces the algebras $A_0$,...,$A_3$
over the real field so that $A=A_0i_0\oplus ... \oplus A_3 i_3$. The
algebra $A$ is the left and right module over ${\bf H}$,
consequently, $i_kA=A$ for each $k$, since $|i_kz| =|z|$ for each
$z\in {\bf H}$. Therefore, the algebras $A_0,...,A_3$ over the real
field are pairwise isomorphic. Moreover, $Z(A)= \emptyset $ implies
that $Z(A_0)= \emptyset $, since $f(z)\ne 0$ means that one of the
components $f_k(z)$ is non-zero, while $A_0$ and $A_k$ are
isomorphic, where $f(z)=\sum_k f_k(z)i_k$ with real-valued
components $f_k$, $f\in A$, $f_k\in A_k$ for each $k$. Since $A_0$
is dense in $C(U,{\bf R})$, then $A$ is dense in $C(U,{\bf H})$.
\par {\bf 5. Corollary.} {\it If $U$ is a canonical closed bounded subset
in ${\bf H}^k$, then the family ${\cal H}(U,{\bf H}^m)$ of all ${\bf
H}$-differentiable functions $f: U\to {\bf H}^m$ is dense in
$C(U,{\bf H}^m)$, where $k, m \in \bf N$.}
\par {\bf Proof.} The set $U$ is closed and bounded in ${\bf H}^k$,
as the $\bf R$-vector topological space ${\bf H}^k$ is locally
compact, since $k \in {\bf N}$, consequently, $U$ is compact.
Therefore, the statement of this corollary follows from Theorem 4,
since ${\cal H}(U,{\bf H}^m)$ is the subalgebra in $C(U,{\bf H}^m)$
and $Z({\cal H}(U,{\bf H}^m)) =\emptyset $.

\par {\bf 6. Theorem.} {\it Suppose that a canonical closed domain $U$ in the
quaternion skew field ${\bf H}$ is compact. Then the set ${\cal
P}(U,{\bf H})$ of all polynomials $P_n: U\to {\bf H}$ is dense in
$C(U,{\bf H})$.}
\par {\bf Proof.} To prove this theorem we use the preceding theorem.
Particularly, we take $T=U$.
\par To rewrite a function from the real variables $z_j$ in the
$z$-representation or vice versa the following identities are used:
\par $(1)$ $z_j=(-zi_j+ i_j2^{-1} \{ -z
+\sum_{k=1}^3i_k(zi_k^*) \} )/2$ \\ for each $j=1, 2, 3$,
\par $(2)$ $z_0=(z+ 2^{-1} \{ -z +
\sum_{k=1}^3i_k(zi_k^*) \} )/2$, \\
where $z$ is a quaternion number decomposed as
\par $(3)$ $z=z_0i_0+...+z_3i_3\in {\bf H}$ \\ with $z_j\in \bf R$ for each $j$,
$i_k^* = {\tilde i}_k = - i_k$ for each $k>0$, $i_0=1$, since
$i_j(i_ji_k)=-i_k$ and $(i_ki_j)i_j=-i_k$ for each $j>0$, also
$i_ji_k=-i_ki_j$ for each $j\ne k$ with $j>0$ and $k>0$, while
$i_k(i_0i_k^*)=1$ for each $k$. Formulas $(1-3)$ define the
real-linear projection operators $\pi _j: {\bf H}\to \bf R$ so that
\par $(4)$ $\pi _j(z)=z_j$ \\ for each quaternion number $z\in
{\bf H}$ and every $j=0,1, 2, 3$.
\par Let $f\in C(U,{\bf H})$. The canonical closed domain $U$ in the
quaternion skew field ${\bf H}$ is compact. It has the real shadow
$V$ in ${\bf R}^4$. Each function $f: U\to {\bf H}$ can be written
in the form
\par $$(5) \quad f(z)=\sum_{j=0}^3 f_j(z) i_j ,$$
where $f_j: U\to {\bf R}$, $~i_j$ is the standard generator of the
quaternion skew field ${\bf H}$ for each $j=0, 1, 2, 3$. On the real
shadow $V$ to each $f_j$ a function $g_j$ of real variables $z_0,
z_1,...,z_3$ corresponds due to equalities $(1-3)$ above. Thus, if
$f: U\to {\bf H}$ is continuous on $U$, then each $g_j: V\to {\bf
R}$ is continuous on $V$. Vice versa if $g_j: V\to {\bf R}$ is
continuous on $V$ for each $j=0,...,3$, then $f: U\to {\bf H}$ is
continuous on $U$ due to formulas $(1-4)$.
\par The set of all real-valued polynomials in the variable
$z=z_0i_0+...+z_3i_3\in U$ forms the algebra $A_0$ over $\bf R$,
since the sum of polynomials and the product of polynomials from
$A_0$ is again a polynomial. Each polynomial of the form \par $(6)$
$q(z) = \pi _j((z-u)^n)$ \\ belongs to $A_0$, where $z, u\in U$,
$~n\in \bf N$, $~u$ is a marked quaternion parameter for $q$. Since
\par
$(7)$
$(z-u)^n==z^n-z^{n-1}u-(z^{n-2}u)z...+(-1)^{n-1}(...((zu)u)...)u+(-u)^n$,
\\ then Formulas $(1-3)$ imply that $q(z)$ is the polynomial with
${\bf H}$ coefficients in the variable $z$, since $\bf H$ is
associative. If $x$ and $y$ are two distinct points in $U$, then
there exists $k\in \{ 0, 1, 2, 3 \} $ such that
$x_k\ne y_k$. Then the function \par $(8)$ $g(z) := (z_k-x_k)^n$, \\
where $n\ge 1$ is a natural number, separates points $x$ and $y$. In
view of Formulas $(1-3)$ this function $g$ expresses as the
real-valued polynomial $Q_{n,k}(z)$ in the variable $z\in U$ so that
$Q_{n,k}(x)=0$ and $Q_{n,k}(y)\ne 0$. Thus the algebra $A_0$
separates points in $U$ such that $Z(A_0)=\emptyset $. Applying
Theorem 4 one gets the statement of this theorem. In more details
the end of the proof is the following.
\par Then in view of the classical Stone-Weierstrass theorem
(see above) $A_0$ is (everywhere) dense in $C(U,{\bf R})$. On the
other hand, Formula $(5)$ means that $C(U,{\bf H})$ is isomorphic
with $[C(U,{\bf R})i_0]\oplus ... \oplus [C(U,{\bf R})i_3]$. The
same Formula $(5)$ implies that the algebra $A={\cal P}(U,{\bf H})$
of all quaternion valued polynomials on $U$ in the quaternion
variable $z\in U$ is isomorphic with $(A_0i_0)\oplus ... \oplus
(A_0i_3)$, since each $P_n(z)\in {\cal P}(U,{\bf H})$ has
the form: \par $(9)$ $P_n(z)=\sum_{j=0}^3 P_{n,j}(z)i_j$ \\
with $P_{n,j}(z)=\pi _j(P_n(z))$ being real-valued polynomials on
$U$ and with ${\bf H}$ expansion coefficients, because ${\bf
R}i_j\subset {\bf H}$ for each $j=0,...,3$. Therefore, ${\cal
P}(U,{\bf H})$ is (everywhere) dense in $C(U,{\bf H})$. That is, for
each $f\in C(U,{\bf H})$ and for each $\epsilon >0$ a polynomial
$P_n(z)\in {\cal P}(U,{\bf H})$ exists such that
$$\| f - P_n \| _{C(U,{\bf H})} < \epsilon .$$

\section{Dirac modules over quaternions.}
\par {\bf 1. Dirac module.} The quaternion skew field $\bf H$ is associative and non-commutative.
It has the standard basis $\{ i_0, i_1, i_2, i_3 \}$
over the real field $\bf R$ such that $i_0=1$, $i_1i_2=i_3$, $i_1^2=i_2^2=i_3^2=-1$, $i_0i_j=i_j$
for each $j$, while
$i_ji_k=-i_ki_j$ for each $j\ne k\ge 1$. We consider the Dirac operator
$\sigma : C^1(\Omega ,{\bf H})\to C(\Omega ,{\bf H})$, where
$\sigma  f(z) = \sum _{j=0}^3 (\partial f(z)/\partial z_j) i_j$
for a differentiable function $f(z)$, $z=\sum _{j=0}^3 z_j i_j$, $z\in {\bf H}$,
$z_j\in {\bf R}$, $z^*=z_0i_0-z_1i_1-z_2i_2-z_3i_3$ denotes the conjugated quaternion
$z$.
\par More generally we consider an associative Clifford algebra ${\cal X}$ of dimension $l$ over ${\bf H}$
and the kernel of the Dirac operator $ker (\sigma )$ in $C^1(\Omega
,{\cal Y})$ or its submodule, when it is indicated. Where $\Omega $
is a domain in ${\cal X}$, while ${ \cal Y } $ denotes a finite
dimensional (two sided) module over $\cal X$ with basis
$q_1,...,q_m$,
$$(1)\quad \sigma f (z) = \sum_{k=1}^l \sum _{j=0}^3 (\partial f(z)/\partial z_{j,k}) e_ki_j,$$  $e_k$
are basic elements in ${\cal X}$ over ${\bf H}$, such that $e_me_k^* = \delta _{m,k}e_1$, $e_m^* = (-1)^{p(m)} e_m$ with
$p(1)=2$, $p(m)\in \{ 1, 2 \} $ for each $2\le m\le n$, $e_ki_j=i_je_k$ and $e_1i_j=i_j$ for each $j$ and $k$,
$$(2)\quad z= \sum_{k=1}^l \sum _{j=0}^3 z_{j,k} e_ki_j$$ with $z_{j,k}\in {\bf R}$ for each $j, k$, whilst $z\in {\cal X}$,
$l$ is a natural number, $$(3)\quad z^*= \sum_{k=1}^l \sum _{j=0}^3
z_{j,k} e_k^*i_j^* .$$ As usually $C^m(\Omega ,{\cal Y})$ and
$C(\Omega ,{\cal Y})$ stand for spaces of $m$ times continuously
differentiable (by all real variables $z_{j,k}$) and continuous
respectively functions on a domain $\Omega $ in ${\cal X}$ with
values in ${\cal Y}$. If $\Omega $ is a real $C^m$-manifold embedded
into $\cal{X}$, then in the standard way using charts of an atlas
$At (\Omega )$ spaces $C^m(\Omega ,{\cal Y})$ and partial
derivatives $D^{\alpha }f$ in local coordinates are defined. For the
unit sphere $S^{n-1}$ the traditional atlas consisting of two charts
is considered (see, for example, \cite{hirschb}).
\par We denote by $K(B,{\cal Y})$ the space of all $C^1$ functions $f$ on $B$ with values in $
{\cal Y}$ such that $\sigma f(z)=0$ for each $z\in B$ and a function
$f(z)$ has a bounded continuous extension on the closure ${\bar
B}=\{  z: z\in {\cal X}, |z|\le 1 \} $ of the open unit ball $B$ in
${\cal X}$. The family $K(B,{\cal Y})$ is considered relative to the
$C$ norm
$$(4)\quad \| f  \|_{C(B,{\cal Y})} := \sup \{ |f(z)|: z\in  {\bar B} \} ,$$
where $ |z|^2 =\sum_{k=1}^m |\mbox{}_kz |^2$ for each
$z=\mbox{}_1zq_1+...+\mbox{}_mzq_m$, $z\in {\cal Y}$, where
$~\mbox{}_kz\in {\cal X}$ for each $k$, whilst $|y|^2 = \sum_{k=1}^l
|y_k|^2$ for any $y=y_1e_1+...+y_le_l\in {\cal X}$ with $y_k\in
{\cal A}_r$ for each $k$.

\par {\bf 2. Proposition.} {\cal The family $K(B,{\cal Y})$ has the structure of the left $\bf H$ module.
If $f\in K(B,{\cal Y})$, then $f$ is harmonic on $B$, that is $\Delta f(z)=0$ for each $z\in B$, where
$\Delta $ denotes the Laplace operator.}
\par {\bf Proof.} From Formula 1$(1)$ it follows that
\par $(1)\quad \sigma ^* \sigma g(z)  = \sigma \sigma ^* g(z)= \Delta g(z) $ for each $g\in C^2(B,{\cal Y})$, where
\par $(2)\quad \Delta f(z) = \sum_{k=1}^l \sum _{j=0}^3 (\partial ^2f(z)/\partial z_{j,k}^2) $ is the Laplace operator
of $4l$ real variables $z_{j,k}$, whilst
$$(3)\quad \sigma ^* f (z) = \sum_{k=1}^l \sum _{j=0}^3 (\partial f(z)/\partial z_{j,k}) e_k^*i_j^*.$$
Then from $\sigma f(z)=0$ for each $z\in B$ it follows that $\sigma ^* \sigma f(z) =0$ on $B$, consequently, each
$f\in K(B,{\cal Y})$ is harmonic, i.e. $\Delta f(z)=0$ on $B$. On the other hand, for each $f\in K(B,{\cal Y})$ and $b\in {\bf H}$
we get that $$(4)\quad \sigma bf(z) =  b \sum_{k=1}^l \sum _{j=0}^3 (\partial f(z)/\partial z_{j,k}) e_ki_j = b[\sigma f(z)],$$
since the quaternion skew field is associative and $\cal Y$ is the left $\bf H$ module. Thus the condition $f\in K(B,{\cal Y})$
implies $bf\in K(B,{\cal Y})$ for each quaternion $b$. Then for every $f, g \in K(B,{\cal Y})$ and $b, c \in {\bf H}$
we infer that $\sigma (bf+cg) = b [\sigma f] + c  [\sigma g]=0$ and $\sigma (bcf) = \sigma (b(cf)) = b\sigma (cf) =
(bc)\sigma f=0$ on $B$. Thus $ K(B,{\cal Y})$ is the left module over the quaternion skew field.
\par {\bf 3. Corollary.} {\it The left modules $K(B,{\cal Y})$ and $K(B,{\cal X})q_1\oplus ...\oplus K(B,{\cal X})q_m$
are isomorphic, where $q_1,...,q_m$ is a basis of $\cal Y$ over
$\cal X$.}

\par {\bf 4. Remark.} The quaternion skew field is also the particular case of the Clifford algebra.
Functions with values in the Clifford algebra $\cal X$ which satisfy
the equation $\sigma f=0$ are called Clifford analytic, where
$\sigma $ is the Dirac operator over $\cal X$ (see also Chapter 2 \S
3 in \cite{gilmurr}). The non trivial Clifford algebras different
from $\bf R$, $\bf C$ and $\bf H$ have divisors of zero. Therefore,
in Section 1 the case was considered, when $K(B,{\cal Y})$ is the
left module over the quaternion skew field (see Proposition 2). But
due to the analog of the Riemann mapping theorem (see \S 2.47 in
\cite{ludnfjms08} and \S 2.1.5.7 in \cite{ludjms07spg}) over
quaternions for a wide class of domains $\Omega $ it is sufficient
to consider the open unit ball $B = \{  z: z\in {\bf H}^n, |z|< 1 \}
$. Generally one can consider a domain $\Omega $ which is $C^{\omega
}$ diffeomorphic with the unit ball $B= \{  z: z\in {\bf H}^n, |z|<
1 \} $ in ${\bf H}^n$, $1\le n\in {\bf N}$, where $C^{\omega }$
denotes the class of all locally analytic functions from $\Omega $
into $B$.
\par It can be lightly seen that there are $f, g \in K(B,{\bf H})$ the product of which is not
in $K(B,{\bf H})$, because
the quaternion skew field is non commutative and $[(\partial f(z)/\partial z_{j,k})i_j]g$ generally
can be not equal to
$[(\partial f(z)/\partial z_{j,k})g]i_j$. Thus apart from the complex case $K(B,{\cal Y})$
is not an algebra even for ${\cal Y}={\bf H}$.
\par Due to Proposition 2 and Corollary 3 it is sufficient to prove that $K(B,{\cal X})$
has a Schauder basis, then the left module $K(B,{\cal Y})$ would
have it as well.

\par {\bf 5. Theorem.} {\it Let ${\Omega }$ be a connected domain in ${\bf H}^n$. Let also
$u: \Omega \to {\bf H}$ be a harmonic function on $\Omega $. If
$|u|$ has a maximum in $\Omega $, then $u$ is constant.}
\par {\bf Proof.} Suppose that $|u|$ has a maximum value $q$ at a point $x\in \Omega $.
Then $u(x)^* u(x)=|u(x)|^2$ and hence there exists $b\in {\bf H}$
such that $bu(x)=q$ and $|b|=1$. Evidently the real part $Re (bu)$
is the harmonic function on $\Omega $ and attains its maximum value
at $x$ according to Theorem 1.4 \cite{axlerb}. Then the equalities
$|bu(z)|=|b| |u(z)| =|u(z)|$ and the inequality $|u(z)|\le q$ imply
that $Im (bu) := bu - Re (bu) =0$ on $\Omega $. Thus $bu$ and hence
$u$ is constant on $\Omega $, since the quaternion skew field ${\bf
H}$ has not divisors of zero.

\par {\bf 6. Note.}
For Clifford analytic functions the analog of the latter theorem is
contained in \S \S 3.28 and 3.30 \cite{gilmurr}.
\par Let $P(x,y)=(1-|x|^2)/|x-y|^{n}$ be the Poisson kernel on the $n$-dimensional over the real field
$\bf R$ unit ball $\cal B$ in the Euclidean space ${\bf R}^n$ with
$n>2$, where $x \ne y\in \bar{\cal B}$. Then for the unit ball $B$
in ${\bf H}^l$, where $l\in {\bf N}$, using the real shadow ${\cal
B}$ of it one gets the Poisson kernel on $B$ with $n=4l$: \par $(1)$
$P(z,\xi ) = (1-|z|^2)/|z-\xi |^{n}$ for each $z\ne \xi \in \bar B$.
We put $$S^{n-1}_+ := \{ z\in S^{n-1}: z_{0,1}> 0 \} \mbox{ and
}S^{n-1}_- := \{ z\in S^{n-1}: z_{0,1}< 0 \} ,$$ where $S^{n-1} =
\partial B$ denotes the unite sphere. We consider the function $U(z)
:= {\hat P}(\chi _{S^{n-1}_+} - \chi _{S^{n-1}_-}) (z)$, where $\chi
_A $ denotes the characteristic function of subset $A$, $\chi _A
(x)=1$ for each $x\in A$, while $\chi _A(x)=0$ for any $x\notin A$,
$\hat P$ is the integral operator:
$$(2)\quad {\hat P} (g)(z) := \int_{S^{n-1}} g(\xi )P(z,\xi )\psi (d\xi ),$$
where $\psi $ denotes the normalized Riemann volume element
(Borel measure) on $S^{n-1}$, $\psi (S^{n-1})=1$.

\par {\bf 7. Theorem.} {\it Let $f$ be a harmonic function $f: B \to {\bf H}$ with $f(0)=0$, where
$B$ is an open unit ball in ${\bf H}^l$, $l\in \bf N$. Then \par
$(1)$ $|f(z)|\le U(|z|i_0e_1)$ \\ for each $z\in B$. Moreover, if
the equality holds in Formula $(1)$ for some nonzero $z\in B$, $z\ne
0$, then there exists a constant $c\in {\bf H}$ of modulus $1$ and
an ${\bf R}$-linear transformation $T$ of ${\bf H}^l$ so that it
induces an orthogonal operator $[T]$ on the real shadow ${\bf R}^n$
of ${\bf H}^l$ such that \par $(2)$ $f=cU\circ T$ on $B$.}
\par {\bf Proof.} A function $f$ can be written in the form $f(z)=\sum_{j=0}^3 f_j(z)i_j$,
where $f_j$ are real-valued functions. Therefore the function $f$ is
harmonic, $\Delta f=0$, if and only if each function $f_j$ is
harmonic. Then we take a constant $b\in {\bf H}$ of modulus $1$ for
a given $z\in B$ such that $|f(z)|=bf(z)$, since ${\bf H}$ is the
associative skew field. The Laplace operator $\Delta $ is invariant
under each transformation $Q$ such that the corresponding operator
$[Q]$ on the real shadow ${\bf R}^n$ is orthogonal, consequently,
the function $f\circ Q^{-1}$ is harmonic as well. Moreover, $\Delta
(bf)=b\Delta f$, since the Laplace operator is the partial
differential operator in real variables with real coefficients,
while the real field is the center of the quaternion skew field
${\bf H}$, hence $bf$ is harmonic, when $f$ is harmonic. One can
take in particular $Q$ such that $[Q]$ is orthogonal and
$Q|z|i_0e_1=z$. In view of Theorem 6.16 \cite{axlerb} applied to the
real part $Re ( bf\circ Q^{-1})$ one gets inequality $(1)$ on $B$.
\par On the other hand, the equality in Formula $(1)$ implies that
$f(z) = b^*U(Q|z|i_0e_1)$. But $|bf(x)|\le U(x)$ for each $x\in B$.
For the real valued function $Re (bf)$ the equality at $z\ne 0$ in
$B$ implies that
\par $(3)$ $Re (bf)(x) = U(x)$ for each $x\in B$ \\ due to Theorem 6.16
\cite{axlerb}. Therefore, $Im (bf(x))=0$ for each $x\in B$ by already proved
Inequality $(1)$ and Property $(3)$. Thus Formula $(2)$ is fulfilled with $c=b^*$.

\par {\bf 8. Theorem.} {\it For each function $f\in K(B,{\cal X})$ there exists a continuous function
$h: S^{n-1}\to {\cal X}$ such that
$$(1)\quad f(z)= \int_{S^{n-1}} h(\xi ) \sigma ^*_zP(z,\xi )\psi (d\xi )$$
for each $z\in B$, where $B$ is an open unit ball with center at
zero in ${\cal X}$ (see \S 1), $S^{n-1} = \partial B$.}
\par {\bf Proof.} According to Proposition 2 from the inclusion $f\in K(B,{\cal X})$ it follows that
a function $f$ is harmonic. Moreover, Formulas 2$(1-3)$ imply that
\par $(2)$ $f= \sigma ^*_zg$ \\ for some harmonic function $g$ on
$B$, since $\Delta =\sigma \sigma ^*$ and the operators $\sigma $
and $\sigma ^*$ commute on $C^2$, while each harmonic function on
$B$ is infinite differentiable by Theorem 1.18 \cite{axlerb}.
Indeed, one equation $(2)$ can be written as the system of $4l$
linear partial differential equations with constant real
coefficients in $4l$ real variables and real functions $f_{j,p},
g_{j,p}$:
$$(3)\quad f_{j,p}(z) = \sum_{t,q; ~ m, k; ~ i_mi_k^*=i_j; ~ e_qe^*_t=e_p}
[\partial g_{m,q}(z)/\partial z_{k,t} + (-1)^{\phi (m,k)} \partial g_{k,q}(z)/\partial z_{m,t}],$$
where $j=0, 1, 2, 3$; $p=1,...,l$; natural numbers $\phi (m,k)\in \{ 1, 2 \}$ are such that
\par $(4)$ $i_mi_k = (-1)^{\phi (m,k)}i_ki_m$. \\
Since $f$ is continuous on $\bar B$,
then the function $g$ is also continuous on $\bar B$. In view of Theorem 1.17  \cite{axlerb} the function $g$ can be written in the form
\par $(5)$ $g(z) = {\hat P}(g|_{S^{n-1}})(z)$, \\ where the integral operator is given by Formula 6$(2)$, $g|_{S^{n-1}}$ denotes the restriction of $g$ on $S^{n-1}$, since the Poisson kernel is real, $P(z,\xi )\in {\bf R}$ for each $z\ne \xi \in \bar B$, and hence the operator ${\hat P}$ is $\bf H$ linear. Evidently the integral in Formula $(1)$ converges uniformly on each smaller closed ball ${\bar B}_R$ of radius $0<R<1$ with center at zero in $B$.
By the theorem of analysis about differentiation of integrals
depending on parameters (see \S XVII.2.3 in \cite{zorichb}) the
identity is valid $$(6)\quad \sigma ^*_z{\hat P}(g|_{S^{n-1}})(z)=
\int_{S^{n-1}} g(\xi ) \sigma ^*_zP(z,\xi )\psi (d\xi )$$ for each
$z\in B$. Thus Formulas $(2,5,6)$, 2$(1)$ and 6$(2)$ lead to the
representation $(1)$.

\par {\bf 9. Remark.} Theorem 8 shows that the left module $K(B,{\cal X})$ is infinite dimensional.
\par For each $y, z \in {\cal X}$ (see \S 1) we put
$$(1)\quad <y,z> := \sum_{p=1}^l y_p^*z_p $$ to be the quaternion valued scalar product.
Let $\mu $ be a Lebesgue measure on the real shadow ${\bf R}^n$,
then $L^2(B,{\cal X})$ denotes the space of all $\mu $-measurable
functions $f: B\to {\cal X}$ such that $\| f \|_2 <\infty $, where
$\| f \|_2:= \sqrt{ (f,f)}$, the $\bf H$ valued scalar product is
given by the integral:
$$(2)\quad (f,g) := \frac{\int_B <f(z),g(z)>\mu (dz)}{\mu (B)}.$$

\par {\bf 10. Theorem.} {\it Let $f\in C^1(\bar{B}, {\cal X})$ be a function harmonic on $B$, and let
$g\in C^1(\bar{B},{\bf R})$ be a real valued function, then
$$(1)\quad (\sigma ^*f,\sigma ^*g) = n\int_{S^{n-1}} <\sigma ^*f(y),y^*> g(y) \psi (dy).$$}
\par {\bf Proof.} Definition 9$(1)$ of the scalar product implies that $<y,z>=y^*z$. We consider the spherical coordinates in ${\bf R}^n$ related with the Cartesian coordinates by the formulas:
\par $(1)$ $x_1=r\cos (\theta _1)$,
\par $x_2= r \sin (\theta _1)\cos (\theta _2)$,...,
\par $x_{n-1} = r \sin (\theta _1)\sin (\theta _2)...\sin (\theta _{n-2})\cos (\theta _{n-1})$,
\par $x_{n} = r \sin (\theta _1)\sin (\theta _2)...\sin (\theta _{n-2})\sin (\theta _{n-1})$,
where $r=|x|\ge 0$, $0\le \theta _j\le \pi $ for each $j=1,...,n-2$, while $0\le \theta _{n-1}\le 2\pi $
(see \S 12.1 in \cite{zorichb}), $x=\phi (z)$, $\phi : {\cal X}\to {\bf R}^n$ is the real linear isometry such that,
$x_1=z_{0,1}$, ...,$x_4=z_{3,1}$,...,$x_{n-3}=z_{0,l}$,...,$x_n=z_{3,l}$, $n=4l$, $z\in {\cal X}$, $x\in {\bf R}^n$.
The Jacobian \par $J= J(r,\theta ) = r^{n-1} \sin ^{n-2}(\theta _1)  \sin ^{n-3}(\theta _2)... \sin (\theta _{n-2})$ \\
is positive for $r>0$ and $0<\theta _j<\pi $ for each $j=1,...,n-2$,
where $\theta =(\theta _1,...,\theta _{n-1})$. Then \par $\mu (dz) =
J dr d\theta _1...d \theta _{n-1} = \psi (dy) dr $, where $y=z/r$
for $r>0$. This transformation from Cartesian to spherical
coordinates can be presented as product of the dilation $x\to rx$
and of $n-1$ orthogonal transformations with matrices
$A_1,...,A_{n-1}$ depending on one angle parameter $\theta
_1,...,\theta _{n-1}$ so that $x=(1,0,...,0)rA_1,...,A_{n-1}$, where
$A_k$ is the $n\times n$ matrix with $1$ as diagonal elements
$(j,j)$ for $j\ne k$ and $j\ne k+1$,
$(A_k)_{k,k}=(A_k)_{k+1,k+1}=\cos (\theta _k)$, $(A_k)_{k,k+1} = -
(A_k)_{k+1,k} = \sin (\theta _k)$, others elements of $A_k$ are
zero, where a vector $x$ is written as the one row matrix. That is
one can consider a sequence of $n$ transformations. Using the chain
rule one can express the Dirac operator in spherical coordinates as:
$$(2)\quad \sigma _t f(z(t)) = \sum_{k=1}^n [\partial f(z(t))/\partial t_k]\alpha _k(t) ,$$
where $(t_1,...,t_n)=(r,\theta _1,...,\theta _{n-1})$, functions
$\alpha _k(t)$ have values in ${\cal X}$,
$$(3)\quad \alpha _{1}(t) = \cos (\theta _1) e_1 i_0 + \sin (\theta _1)\cos (\theta _2)e_1 i_1 +...+ \sin (\theta _1)\sin (\theta _2)...\sin (\theta _{n-2})\sin (\theta _{n-1})e_li_3,$$
$$(4)\quad \alpha _{k}(t)=r^{n-2}J^{-1} \sum_{j=0}^3 \sum_{p=1}^l \beta _{k,j,p}(\theta ) e_pi_j ,$$ where $\beta _{k,j,p}(\theta )$ are definite products
$\beta _{k,j,p}(\theta ) = \prod_{m=1}^{n-1}\sin ^{a(m)}(\theta
_m)\cos ^{b(m)}(\theta _m)$, where $a(m)=a_{k,j,p}(m)$ and
$b(m)=b_{k,j,p}(m)$ are nonnegative integers for each $m=1,...,n-1$.
Thus $J(t)\alpha _{k}(t)$ are infinite differentiable functions for
each $k$. In particular, for $l=1$ one has: \par $(5)$ $\alpha _2=
r^2J^{-1}[-\sin ^3(\theta _1)\sin (\theta _2)i_0+ \sin ^2(\theta
_1)\cos(\theta _1) \sin (\theta _2)\cos(\theta _2)i_1$\par $+\sin
^2(\theta _1)\cos(\theta _1)\sin ^2(\theta _2)\cos(\theta _3)i_2+
\sin ^2(\theta _1)\cos(\theta _1)\sin ^2(\theta _2)\sin (\theta
_3)i_3 ]$
\par $\alpha _3= r^2J^{-1}[-\sin (\theta _1)\sin ^2(\theta _2)i_1+ \sin (\theta _1)\sin (\theta _2)\cos(\theta _2)\cos(\theta _3)i_2$\par $+ \sin (\theta _1)\sin (\theta _2)\cos(\theta _2)\sin (\theta _3)i_3 ]$
\par $\alpha _4= r^2J^{-1}[-\sin (\theta _1)\sin (\theta _3)i_2 + \sin (\theta _1)\cos(\theta _3)i_3 ].$
\par In spherical coordinates the adjoint Dirac operator is:
$$(6)\quad \sigma ^*_t f(z(t)) = \sum_{k=1}^n [\partial f(z(t))/\partial t_k]\alpha ^*_k(t) .$$
The operator $\sigma ^*$ in Cartesian coordinates or correspondingly $J\sigma ^*$ in spherical coordinates
defines a vector field $Y^*$ with coefficients in ${\cal X}$. For the Dirac operator $\sigma $ in Cartesian
coordinates the corresponding vector field
$$Y = \sum_{k=1}^l \sum _{j=0}^3 e_ki_j (\partial /\partial z_{j,k}) $$ has constant Clifford coefficients
$e_ki_j$ of unit norm. Then \par $\sigma ^*[(\sigma ^*f)^*g]=<\sigma
^*f, \sigma ^*g>$, \\ since $\Delta f=0$ and the function $g$ is
real valued, while $\bf R$ is contained in the center of the
Clifford algebra $\cal X$. Under the composition of the mappings
$\phi $ and $x\mapsto t =(r,\theta )$ the images of $\bar{B}$ and
$S^{n-1}$ are $Q$ and $\partial Q$ respectively, where $Q :=
[0,1]\times [0,\pi ]^{n-2}\times [0,2\pi ]$, consequently,
$$ (\sigma ^*f,\sigma ^*g) =\int_B \sigma ^*[(\sigma ^*f)^*(z) g(z)] \mu (dz) =
\int_Q \sigma ^*[(\sigma ^*f)^*(z((r,\theta )) g(z(r,\theta ))] J dr
d\theta _1...d\theta _{n-1}.$$ For an even dimensional unit ball
$B=B_n$ its volume $\mu (B) = V_n(B) = \pi ^{n/2} /(n/2)!$ and the
unnormalized surface area of the unit sphere $S^{n-1}$ in ${\bf
R}^{n-1}$ is $nV_n(B)$, where $n$ is the dimension of $B$ over $\bf
R$, $\mu ([0,1]^n)=1$, (see Appendix A in \cite{axlerb} or
\cite{zorichb}). In view of Stokes' theorem XV.3.5 \cite{zorichb}
applied to the integrals in the latter formula one gets $$ (\sigma
^*f,\sigma ^*g) = \sum_{p=1}^n \int_{\partial Q_p} <\sigma
^*f(z(t)),\gamma ^*(t)> g(z(t)) dt_1\wedge ...\wedge dt_{p-1}\wedge
dt_{p+1}...\wedge dt_n$$
$$ = n\int_{S^{n-1}} <\sigma ^*f(y),y^*>g(y)\psi (dy),$$ since $y\in
S^{n-1}$ is a vector orthogonal to $S^{n-1}$ at $y$ and directed
outwards $B$ and of unit norm, $|y|=1$, whilst $(ab)^*=b^*a^*$ for
each $a, b \in {\cal X}$, where $$d\gamma |_{\partial Q}  = J
\sum_{p=1}^n (-1)^{p+1}\alpha _p \chi _{\partial Q_p}dt_p,\quad
\partial Q_p := \{ t\in \partial Q: t_p=0\mbox{ or } t_p=b_p \} ,$$
where the orientation of $Q$ is consistent with that of $\partial Q$
and the orientation of $B$ is consistent with that of $\partial
B=S^{n-1}$ (see also \S XV.3 \cite{zorichb}), $b_1=1$, $b_p=\pi $
for $p=2,...,n-1$, $b_n=2\pi $.

\par {\bf 11. Note.} Using the identity 10$(1)$ we define the scalar product on $L^2(S^{n-1},{\bf R})\cap C^1(S^{n-1},{\bf R})$:
$$(1)\quad [f,g] := n \int_{S^{n-1}} <\sigma ^*\hat{P}(f)(y),y^*> \hat{P}(g)(y)\psi (dy).$$
Let the operator $T$ on $C^1(S^{n-1},{\cal X})$ be defined by the
formula
\par $(2)$  $Tf:= \sigma ^*\hat{P}[f]$.
\par {\bf 12. Corollary.} {\it The restriction of the operator $T$ (see 11$(2)$) to $C^1(S^{n-1},{\bf R})$
induces the real linear isometry from $V:=(L^2(S^{n-1},{\bf R})\cap
C^1(S^{n-1},{\bf R}))/Y$ into $L^2(\bar{B},{\cal X})$ relative to
the corresponding scalar products $[*,*]$ and $(*,*)$, where
$Y=T^{-1}(0)=ker(T)$.}
\par {\bf Proof.} The calculation of $\sigma ^*_zP(z,w)$ gives
$$(1)\quad \sigma ^*_zP(z,w) = - \frac{2|z-w|^2z^*+n(1-|z|^2)(z-w)^*}{|z-w|^{n+2}}$$
for each $z\ne w\in \bar{B}$. Therefore, for each $f\in
L^2(\bar{B},{\cal X})$ we have $\sigma ^*\hat{P}(f)(z)\in {\cal X}$
for any $z\in B$ and Formula 9$(2)$ leads to the inequality $(f,f)
\ge 0,$ since $<f(z),f(z)>\ge 0$ for each $z\in B$, whilst $\mu
(dz)$ is the nonnegative Lebesgue measure. There is the inclusion
$C^1(S^{n-1},{\cal X})\subset L^2(S^{n-1},{\cal X})$. If $f\in
C^1(S^{n-1},{\bf R})$, then the function $\hat{P}(f)$ is harmonic on
$B$ and continuously differentiable on $\bar B$ due to Theorems 1.14
and 1.17 \cite{axlerb}. \par At the same time, the equality
$(f,f)=0$ is equivalent to $<f(z),f(z)>=0$ for almost all $z\in B$,
that in its turn is equivalent to $f(z)=0$ for almost all $z\in B$
according to 9$(1)$, consequently, $\sqrt{(f,f)}$ is the norm on
$L^2(\bar{B},{\cal X})$. Therefore, $\sqrt{[g,g]}=\sqrt{(Tg,Tg)}$ is
the norm on the real quotient space $V$. \par If $g\in
C^1(S^{n-1},{\bf R})$, then $<\sigma ^*\hat{P}(g)(y),y^*> \in {\cal
X}$ for each $y\in S^{n-1}$ as follows from Formulas $(1)$,
1$(1,3)$, 2$(3)$ and 9$(1)$. In view of Equalities $(1)$, 9$(2)$ and
10$(1)$ the inclusion $[f,g]\in {\bf R}$ is valid for each $f, g\in
C^1(S^{n-1},{\bf R})$ due to the polarization identity
$$[f,g]= \frac{[f+g,f+g] - [f-g,f-g]}{4},\mbox{ since }$$
$$(\sigma ^*_zP(z,w))(\sigma _zP(z,\xi )) +(\sigma ^*_zP(z,\xi
))(\sigma _zP(z,w )) = \Delta _z [P(z,w)P(z,\xi )]\in {\bf R}$$ and
$(\sigma ^*_zP(z,w))^*= \sigma _zP(z,w)$ for every $w, \xi , z\in
\bar{B}$ with $w\ne z$ and $\xi \ne z$. On the other hand,
$\hat{P}(g)$ is the harmonic function on $B$. This implies that
$[f,f]=0$ if and only it $f\in ker (T)$. Thus $T: V\to
L^2(\bar{B},{\cal X})$ is the isometry, which is linear over the
real field, where $V$ is supplied with the scalar product $[f,g] =
(Tf,Tg)$ and the norm $ \| g \| = \sqrt{[g,g]}$.

\par {\bf 13. Theorem.} {\it On the Banach space $C^1(S^{n-1},{\bf R})$ the integral
$$(1)\quad (\sigma ^*_zP(z,\xi ),\sigma ^*_zP(z,w)) := \int_B <\sigma ^*_zP(z,\xi ),\sigma ^*_zP(z,w)>\mu (dz)$$
defines the generalized function $(\sigma ^*_z\hat{P}(\delta _{z,\xi
})g(\xi ),\sigma ^*_z\hat{P}(\delta _{z,w})f(w))$ (continuous ${\bf
R}$-bilinear functional with values in ${\cal X}$) by $f$ and $g$,
where $f, g \in C^1(S^{n-1},{\bf R})$.}
\par {\bf Proof.} The delta function $\delta _{\xi ,y}$ on the sphere $S^{n-1}$ is characterized by the formula
$$(2)\quad \int_{S^{n-1}}\delta _{\xi ,y} g(y)\psi (dy) = \int_{S^{n-1}}\delta _{y,\xi } g(y)\psi (dy) = g(\xi )$$
for each $g\in C(S^{n-1},{\bf R})$ and $\xi \in S^{n-1}$.
 Take any delta sequence of continuous nonnegative functions $g_m: S^{n-1}\to {\bf R}$ such that
$$\int_{S^{n-1}}g_m(y)\psi (dy) =1$$ and there exists $0<\epsilon <1$ for which
$$\int_{y\in S^{n-1}: |y-w| < \epsilon /m}g_m(y)\psi (dy) \ge 1-1/m$$ for each $m\in {\bf N}$, that is
the limit exists $\lim_m g_m(y) = \delta _w(y)$ relative to the
weak* topology in the topological dual space $C^*(S^{n-1},{\bf R})$
(see \S 6.4 \cite{narib}), where
$$\int_{S^{n-1}} \delta _w(y) g(y)\psi (dy) = g(w)$$ for each $g\in C(S^{n-1},{\bf R})$,
$(C,C^*)$ is the dual pair. In particular this is for $g(y)=\sigma
^*_zP(z,y)$, when $z\in B$ and $y\in S^{n-1}$. Since
$\hat{P}[f](y)=f(y)$ for each $f\in C(S^{n-1},{\bf R})$ and $y\in
S^{n-1}$ (see Theorems 1.14 and 1.17 in \cite{axlerb}) taking the
limit one arrives to the equality
\par $(3)$ $\hat{P}[\delta _w]|_{S^{n-1}}=\lim_m\hat{P} [g_m]|_{S^{n-1}}=\lim_m g_m|_{S^{n-1}}=\delta _w|_{S^{n-1}}$ \\ in $C^*(S^{n-1},{\bf R})$.
On the other hand, $\int_{S^{n-1}} \sigma ^*_zP(z,y)\delta _w(y)\psi (dy) = \sigma ^*_zP(z,w)$ for each $z\in B$ and $w\in S^{n-1}$.
\par For the direct product of delta functions and every $f, g \in C(S^{n-1},{\bf R})$
the equality $(\delta _{\xi ,y}\times \delta _{w,y})(f(\xi )\times g(w))= f(y)g(y)$ is valid.
Therefore,
by Fubini's theorem $$\int_{S^{n-1}} f(\xi )\int_{S^{n-1}} (\delta _{\xi ,y}\times \delta _{w,y})
\psi (dy) \psi (d\xi )=\int_{S^{n-1}\times S^{n-1}} f(\xi )(\delta _{\xi ,y}\times \delta _{w,y})
\psi (dy) \psi (d\xi )$$
$$=f(w)= \int_{S^{n-1}} \delta _{\xi ,w}f(\xi )\psi (d\xi ),$$ consequently,
$\int_{S^{n-1}} (\delta _{\xi ,y}\times \delta _{w,y})\psi (dy) = \delta _{\xi ,w}$.
Therefore, from Fubini's theorem, Theorem 10 and Equality $(3)$ we infer for each
$f \in C^1(S^{n-1},{\bf R})$
that $$\int_{S^{n-1}}(\sigma ^*_zP(z,\xi ),\sigma ^*_zP(z,w))f(w)\psi (dw) :=$$
$$ \int_{S^{n-1}}\int_B <\int_{S^{n-1}} (\sigma ^*_zP(z,y))\delta _{\xi ,y}\psi (dy),
\int_{S^{n-1}}\sigma ^*_zP(z,q)\delta _{w,q}\psi (dq)>f(w)\mu (dz)\psi (dw) $$
$$= n \int_{S^{n-1}} \int_{S^{n-1}}\int_{S^{n-1}} \int_{S^{n-1}}
<(\sigma ^*_zP(z,y)),z^*>\delta _{\xi ,y}\times \delta _{w,q}f(w)P(z,q)\psi (dy)\psi (dq)
\psi (dw) \psi (dz)$$
$$= n \int_{S^{n-1}} \int_{S^{n-1}}\int_{S^{n-1}}
<(\sigma ^*_zP(z,y)),z^*>\delta _{\xi ,y}\times \delta _{w,z}f(w)\psi (dy)\psi (dw) \psi (dz)$$
$$= n \int_{S^{n-1}} \int_{S^{n-1}}
<(\sigma ^*_wP(w,y)),w^*>f(w)\delta _{\xi ,y}\psi (dy)\psi (dw) $$
$$= n \int_{S^{n-1}} <(\sigma ^*_wP(w,\xi )),w^*>\hat{P}(f)(w)\psi (dw) $$
$$= \int_{B} <(\sigma ^*_wP(w,\xi )),(\sigma ^*_w\hat{P}(f)(w))>\mu (dw) =
(\sigma ^*_z\hat{P}(\delta _{z,\xi }),\sigma ^*_z\hat{P}(f)(z)).$$
Finally we get $(\sigma ^*_z\hat{P}(\delta _{z,\xi }),\sigma
^*_z\hat{P}(f)(z))g(\xi )= (\sigma ^*_z\hat{P}(\delta _{z,\xi
})g(\xi ),\sigma ^*_z\hat{P}(\delta _{z,w})f(w))$. The continuity of
this real bilinear functional follows from the equalities
$$(\sigma ^*_z\hat{P}(\delta _{z,\xi })g(\xi ),\sigma ^*_z
\hat{P}(\delta _{z,w})f(w))= [f,g] := n \int_{S^{n-1}} <\sigma ^*\hat{P}(f)(y),y^*> \hat{P}(g)(y)\psi (dy)$$ and the estimate
$$|[f,g]| \le n \| f \|_{C^1(S^{n-1},{\bf R})}  \| g \|_{C^1(S^{n-1},{\bf R})},$$ since $\sigma ^*\hat{P}(f)(y)|_{S^{n-1}} = \sigma ^*f(y)|_{S^{n-1}}=
\hat{P}(\sigma ^*f(y)|_{S^{n-1}})|_{S^{n-1}}$.

\par {\bf 14. Proposition.} {\it The operators $\hat{P}: C(S^{n-1},{\cal X})\to K(B,{\cal X})$ and $\sigma ^*\hat{P}: C^1(S^{n-1},{\cal X})\to K(B,{\cal X})$ are continuous, $\hat{P}$ is left and right $\bf H$ linear, while $\sigma ^*\hat{P}$ is left $\bf H$ linear.}
\par {\bf Proof.} The integral operator $\hat{P}$ has the real integral kernel $P$ and the Borel measure $\psi $ on $S^{n-1}$ in the integral is nonnegative, hence $\hat{P}$ is left and right $\bf H$ linear. In \S 2 it was proved, that the operator $\sigma ^*$ is left $\bf H$ linear, consequently,
the composite operator $\sigma ^*\hat{P}$ is left $\bf H$ linear.
Since $f(z)\mapsto f(z)-f(0)$ is the continuous mapping of
$C(S^{n-1},{\bf H})$ into itself, $\hat{P}(f)$ is harmonic on $B$
and $\hat{P}(f)|_{S^{n-1}}=f|_{S^{n-1}}$ for each $f\in
C(S^{n-1},{\bf H})$, then from Theorem 7 it follows that the mapping
$\hat{P}: C(S^{n-1},{\bf H})\to K(B,{\bf H})$ is continuous.
\par On the other hand, $K(B,{\cal X})=\bigoplus_{p=1}^l  K(B,{\bf H})e_p$, hence the operator $\hat{P}$ from $C(S^{n-1},{\cal X})$ into $K(B,{\cal X})$ is also continuous. Analogously the mapping $\hat{P}: C^1(S^{n-1},{\cal X})\to C^1(\bar{B},{\cal X})$ is continuous, since each harmonic function
is infinite differentiable. The operator $\sigma ^*:
C^1(\bar{B},{\cal X})\to C(\bar{B},{\cal X})$ is continuous,
consequently, the operator $\sigma ^*\hat{P}$ from
$C^1(S^{n-1},{\cal X})$ into $K(B,{\cal X})$ is continuous as well.

\par {\bf 15. Theorem.} {\it The left module $K(B,{\cal X})$ is Banach and
has a Schauder basis over the quaternion skew field $\bf H$.}
\par {\bf Proof.} By the left $\bf H$ span $l-span_{\bf H} \Psi $ of a set $\Psi $ in a
left module over $\bf H$ we mean all finite sums $a_1f_1+...+a_lf_l$ with
$a_1,...,a_l\in {\bf H}$ and $f_1,...,f_l\in \Psi $, where $l\in {\bf N}$.
At the same time the operators $\sigma $ and $\sigma ^*$ map
$C^1({\bf R}^n,{\cal X})$ into $C({\bf R}^n,{\cal X})$ and hence
$C^1(S^{n-1},{\cal X})$ into $C(S^{n-1},{\cal X})$. According to
Tietze-Urysohn's theorem (see \S 2.1.8 in \cite{eng}) each
continuous function $u:  S^{n-1}\to {\bf R}$ has a continuous
extension $u: \bar{B}\to {\bf R}$, consequently, each continuous
function $f: S^{n-1}\to {\cal X}$ has a continuous extension $f:
\bar{B}\to {\cal X}$, since $\bar B$ is the normal topological space
and $S^{n-1}$ is closed in it. \par In view of Theorem 8 and
Formulas 10$(2-4)$ we have, that
$\hat{P}(f|_{S^{n-1}})|_{S^{n-1}}=f|_{S^{n-1}}$ and
$\hat{P}(f|_{S^{n-1}})$ is harmonic on $B$ so that if $g\in
C^1(S^{n-1},{\cal X})$, then \par $(1)$ $\sigma
^*\hat{P}(g|_{S^{n-1}})|_{S^{n-1}}=\sigma ^*g|_{S^{n-1}}$. If $f\in
C(\bar{B},{\cal X})$, then the linear system of partial differential
equations with real coefficients 8$(3)$ has a solution $g\in
C^1(\bar{B},{\cal X}))$ (see Theorem 2.4.1 \cite{ludcmft12} and
references therein), hence $\sigma ^*(C^1(S^{n-1},{\cal X})) =
C(S^{n-1},{\cal X}).$ This implies that $[l-span_{\bf H} \sigma
^*(C^1(S^{n-1},{\bf R}))]\cap C(S^{n-1},{\bf H})$ is dense in
$C(S^{n-1},{\bf H})$ and hence $\bigoplus_{p=1}^l l-span_{\bf H}
[\sigma ^*(C^1(S^{n-1},{\bf R})]e_p$ is dense in $C(S^{n-1},{\cal
X})$. At the same time we have from Formula 2$(3)$, that the kernel
of the restriction to $C^1(\bar{B},{\bf R})$ of the operator $\sigma
^*$ consists of all constant functions, that is $ker (\sigma
^*|_{C^1(\bar{B},{\bf R})} ) = {\bf R}$.
\par Then we take a Schauder basis $\{ f_m: ~ m\in {\bf N} \} $ in $C(S^{n-1},{\bf H})$
(see Theorem 2.4 above), consequently, its left $\bf H$ span is
dense in $C(S^{n-1},{\bf H})$. Therefore, $\{ f_me_k: ~ m\in {\bf
N}, k=1,...,l \} $ is the Schauder basis in the left $\bf H$ module
$C(S^{n-1},{\cal X}).$ For each $f_m$ then we choose a particular
solution $g_m\in C^1(S^{n-1},{\cal X})$ of the equation $\sigma
^*g_m=f_m$ restricted to $S^{n-1}$ as above (see \S 8 also). We
modify this basis $ \{ f_me_k: m, k \} $ in such manner that each
$g_m$ is real valued, which is possible, since $C^1(S^{n-1},{\cal
X})=\bigoplus_{p=1}^l[C^1(S^{n-1},{\bf R})i_0\oplus ...\oplus
C^1(S^{n-1},{\bf R})i_3]e_p.$ Each function $\hat{P}(f_m)$ and
$\hat{P}(g_m)$ is harmonic on $B$, moreover, each $\hat{P}(g_m)$ is
real valued.
\par From Formulas 1$(2)$, 9$(1,2)$, 10$(1)$, 11$(1)$ and Corollary 12 it follows that
$[g_m,g_p] \in {\bf H}$ for any $m, p \in {\bf N}$. Each quaternion
equation $ax=b$ or $xa=b$ has a solution $x\in {\bf H}$, when $a\ne
0$, where $a, b\in {\bf H}$. Then using Schmidt's orthogonalization
and normalization procedures relative to the scalar product $[*,*]$
applied to $\{ g_m: m \} $ we get functions
$u_m=\sum_{k=1}^ma_{m,k}f_k$ and $v_m=\sum_{k=1}^ma_{m,k}g_k$ such
that $[v_m,v_p] = \delta _{m,p}$ for each $m, p \in {\bf N}$, where
$a_{m,k}\in {\bf H}$ are quaternion constants, since each $g_k$ is
real valued, while $\delta _{m,p}$ is Kroneker's delta-symbol,
$\delta _{m,p}=0$ when $m\ne p$ whilst $\delta _{m,m}=0$ for every
$m, p\in {\bf N}$. Certainly, these functions are related by the
equation $u_m=\sigma ^*v_m$ for each natural number $m$. Thus each
function $v_m$ is $\bf H$ valued.
\par The operator $\hat P$ is continuous from $C(S^{n-1},{\cal X})$ into $K(B,{\cal X})$.
At the same time a solution of the equation $\sigma ^*g=h$ depends
continuously on $h\in C(\bar{B},{\cal X})$ as it is known from the
theory of systems of linear partial differential equations with
constant coefficients. This means that the anti-derivation operator
$\Upsilon _{\sigma ^*}$ is continuous from $C(\bar{B},{\cal X})$
into $C^1(\bar{B},{\cal X})$ (see Theorem 2.4.1 in \cite{ludcmft12}
and references therein). Therefore, there exists a continuous
anti-derivation operator denoted by $\Upsilon _{\sigma
^*}|_{S^{n-1}}$ from $C(S^{n-1},{\cal X})$ into $C^1(S^{n-1},{\cal
X})$. \par Let $Q$ be the left Banach module over $\bf H$ which let
be the closure in $C^1(S^{n-1},{\cal X})$ of $l-span_{\bf H} \{
v_me_k: ~ m\in {\bf N}, k=1,...,l \} $. Therefore, this Banach
module $Q$ is contained in $L^2(S^{n-1},{\cal X})$. By the
construction above we infer that $\sigma ^*Q=C(S^{n-1},{\cal X})$.
The continuity of the anti-derivation operator $\Upsilon _{\sigma
^*}|_{S^{n-1}}$ implies that $\{ v_me_k: ~ m\in {\bf N}, k=1,...,l
\} $ is the Franklin system in $Q$ relative to the scalar product
$[*,*]$, since
\par $(2)$ $\sigma ^*\hat{P} (\Upsilon _{\sigma ^*}|_{S^{n-1}}f)|_{S^{n-1}} = f= \hat{P}(f)|_{S^{n-1}}$
\\ for each continuous function $f: S^{n-1}\to {\cal X}$, also
$Q\subset (L^2(S^{n-1},{\cal X})\cap C^1(S^{n-1},{\cal X}))$ and
$[v_m,v_p] = \delta _{m,p}$ for each $m, p \in {\bf N}$.
\par Suppose that $h_n$ is a fundamental sequence in $K(B,{\cal X})$, then
$h_n|_{S^{n-1}}$ is a fundamental sequence in $C(S^{n-1},{\cal X})$,
where $n\in {\bf N}$. Therefore, the limit $\lim_{n\to \infty }
h_n|_{S^{n-1}} = y$ exists in $C(S^{n-1},{\cal X})$. Then from the
continuity of the operator $\sigma ^*\hat{P} \Upsilon _{\sigma
^*}|_{S^{n-1}}$, Identities $(2)$ and Theorem 8 we infer that the
limit $$\lim_{n\to \infty } h_n = h= \sigma ^*\hat{P} (\Upsilon
_{\sigma ^*}|_{S^{n-1}}y)$$ exists in $K(B,{\cal X})$ relative to
the $C(\bar{B},{\cal X})$ norm. Hence $K(B,{\cal X})$ is the Banach
left module over the quaternion skew field.
\par If $h\in K(B,{\cal X})$, then from the embedding
$K(B,{\cal X})\hookrightarrow L^2(\bar{B},\mu ,{\cal X})$ and from
Theorems 10,13 and Corollary 12 it follows that
$$h= \sum_{m=1}^{\infty } \sum_{k=1}^l \beta _{m,k} w_{m,k},$$ where
$w_{m,k} = Tv_me_k$, $$h=\sum_{k=1}^l h_ke_k$$ with $h_k\in K(B,{\bf
H})$ for each $k=1,...,l$, $\beta _{m,k} = (h_k^*,w_{m,k}^*)$. That
is, each expansion coefficient $\beta _{m,k}$ for a function $h$ is
unique. Since
$$(h_k^*,w_{m,k}^*)= \frac{\int_B h_k(z)w_{m,k}^*(z)\mu (dz)}{\mu
(B)},$$ each functional $\beta _{m,k}: K(B,{\bf H})\to {\bf H}$ is
left $\bf H$ linear, where $\beta _{m,k}=\beta _{m,k}(h_k)$, since
the quaternion skew field is associative. From
Cauchy-Bunyakovsky-Schwarz's inequality we infer that
$$|(h_k^*,w_{m,k}^*)|\le \| h_k \|_{L^2(\bar{B},{\bf H})} \| w_{m,k}
\|_{L^2(\bar{B},{\bf H})} \le \| h_k \|_{C(\bar{B},{\bf H})} ,$$
since $\| w_{m,k} \|_{L^2(\bar{B},{\bf H})} =1$ and $\mu (B) <\infty
.$ Therefore, each functional $\beta _{m,k}$ is continuous from
$(K(B,{\bf H}), \| * \|_{C(\bar{B},{\bf H})} )$ into ${\bf H}$ and $
\| \beta _{m,k} \| \le 1$ for every $m$ and $k$. The operator
$\sigma^*\hat{P}\Upsilon_{\sigma^*}|_{S^{n-1}}$ is continuous from
$C(S^{n-1},\cal{X})$ onto $K(B,\cal{X})$ due to Formulas $(1,2)$,
Theorem 8 and the proof above. Thus $\{ w_{m,k}: ~ m\in {\bf N},
k=1,...,l \} $ is the Schauder basis in $K(B,{\cal X})$. Moreover,
it is the Franklin system due to Corollary 12.

\par {\bf 16. Corollary.} {\it The left module $K(B,{\cal Y})$ has Schauder bases.}
\par {\bf Proof.} From the construction of the Schauder basis in $K(B,{\cal X})$
and from Remark 4 we get that there exists a Schauder basis in
$K(B,{\cal Y})$ as well.

\par {\bf 17. Notation.} Let $H^p(\bar{B},\cal{Y})$ denote the
(Hardy) left $\bf H$-module of all measurable functions $f:
\bar{B}\to {\cal Y}$ satisfying the condition:
$$(1)\quad \| f \|_{H^p(\bar{B},{\cal Y})} := [\sup_{0<r\le 1}
r^{1-n} \int_{S ^{n-1}} \| f (ry) \|_{\cal Y} \psi (dy)
]^{1/p}<\infty ,$$ where $1<p<\infty $. As usually $W^p_m(\Omega
,\cal{Y})$ stands for the Sobolev space of all functions $f: \Omega
\to \cal{Y}$ such that their partial derivatives $D^{\alpha }f$ are
measurable for each $|\alpha |\le m$ and
$$(2)\quad \| f \|_{W^p_m(\Omega ,{\cal Y})} := [\sum_{|\alpha |\le m}
\int_{\Omega } \| D^{\alpha } f (y) \|_{\cal Y} \psi (dy)
]^{1/p}<\infty ,$$ where $1<p<\infty $, $\Omega $ is a Riemann $C^m$
manifold, $\psi $ is a Borel measure (that is a volume element) on
$\Omega $. Particularly, one gets $W^p_0(\Omega ,{\cal
Y})=L^p(\Omega ,{\cal Y})$ for $m=0$. Then we denote by $K_p(B,{\cal
Y})$ the space of all functions $f: \bar{B}\to {\cal Y}$ satisfying
the conditions:
\par $(3)$ $f|_B\in C^1(B,{\cal Y})$ and
\par $(4)$ $\sigma f(z)=0$ for each $z\in B$ and \par $(5)$ $f\in
H^p(\bar{B},{\cal Y})$. This left ${\bf H}$-module $K_p(B,{\cal Y})$
is supplied with the norm inherited from $H^p(\bar{B},{\cal Y})$.

\par {\bf 18. Theorem.} {\it If  $1<p<\infty $, then the left module $K_p(B,{\cal X})$
is Banach and has an unconditional basis over the quaternion skew
field $\bf H$.}
\par The {\bf Proof} of this theorem is analogous to that of Theorem
15 and Corollary 16 with the following modifications. \par In view
of Theorem 6.12 \cite{axlerb} the mapping $f\mapsto \hat{P}(f)$ is
the surjective isometry from $L^p(S^{n-1},{\bf R})$ onto
$H^p(\bar{B},{\bf R})$, where $1<p<\infty $. On the other hand, the
restriction of the Dirac operator $\sigma ^*$ to the unit sphere
$S^{n-1}$ maps the Sobolev space $W^p_1(S^{n-1},\cal{X})$ into
$L^p(S^{n-1},\cal{X})$, whilst the anti-derivation operator
$\Upsilon _{\sigma ^*}|_{S^{n-1}}$ restricted to $S^{n-1}$ maps the
Lebesgue space $L^p(S^{n-1},{\cal X})$ into $W^p_1(S^{n-1},{\cal
X})$. Then the linear system of partial differential equations with
real coefficients 8$(3)$ has a solution $g\in W^p_1(\bar{B},{\cal
X})$ for each $f\in L^p(\bar{B},{\cal X})$, hence $\sigma
^*(W^p_1(S^{n-1},{\cal X})) = L^p(S^{n-1},{\cal X}).$ Then from
Formula 2$(3)$ we infer, that the kernel of the restriction to
$W^p_1(\bar{B},{\bf R})$ of the operator $\sigma ^*$ consists of all
functions $f$ each partial derivative $\partial f(z)/\partial
z_{m,k}$ of which is zero almost everywhere on $\bar{B}$. From
Lebesgue's theorems (see Theorems 2 and 3 in \S VI.4 \cite{kolfomb})
applied by each real variable $z_{m,k}$ it follows that such $f\in
ker(\sigma ^*)$ is almost everywhere constant on $\bar{B}$, since
$$\int_0^{x_{m,k}} (\partial f(z)/\partial z_{m,k}) dz_{m,k} =
f(z+(x_{m,k}-z_{m,k})e_ki_m)-f(0)$$ for each
$z+(x_{m,k}-z_{m,k})e_ki_m\in \bar{B}$, where $z\in \bar{B}$. Thus
we get that $ker (\sigma ^*|_{W^p_1(\bar{B},{\bf R})} ) = {\bf R}$.
Hence \par $(1)$ $\sigma ^*\hat{P}(g|_{S^{n-1}})|_{S^{n-1}}=\sigma
^*g|_{S^{n-1}}$ almost everywhere on $S^{n-1}$ for each $g\in
W^p_1(S^{n-1},{\cal X})$. Moreover, we deduce that
\par $(2)$ $\sigma ^*\hat{P} (\Upsilon _{\sigma
^*}|_{S^{n-1}}f)|_{S^{n-1}} = f= \hat{P}(f)|_{S^{n-1}}$\\ almost
everywhere on $S^{n-1}$ for each $f\in L^p(S^{n-1},\cal{X})$.
Therefore, $[l-span_{\bf H} \sigma ^*(W^p_1(S^{n-1},{\bf R}))]\cap
L^p(S^{n-1},{\bf H})$ is dense in $L^p(S^{n-1},{\bf H})$ and hence
$\bigoplus_{p=1}^l l-span_{\bf H} [\sigma ^*(W^p_1(S^{n-1},{\bf
R})]e_p$ is dense in $L^p(S^{n-1},{\cal X})$. Then we infer that the
composite operator $\sigma^*\hat{P}\Upsilon_{\sigma^*}|_{S^{n-1}}$
is continuous from $L^p(S^{n-1},\cal{X})$ onto $K_p(B,{\cal X})$ due
to Theorem 6.12 \cite{axlerb}, Formulas 15$(1,2)$ and Theorem 8.
Evidently, there is the continuous embedding of $K_p(B,{\cal X})$
into $L^p(\bar{B},{\cal X})$.
\par Let $h_n$ be a fundamental sequence in $K_p(B,{\cal X})$, then
$h_n|_{S^{n-1}}$ is a fundamental sequence in $L^p(S^{n-1},{\cal
X})$ according to Formula 17$(1)$. Hence the limit $\lim_{n\to
\infty } h_n|_{S^{n-1}} = y$ exists in $L^p(S^{n-1},{\cal X})$.
Using the continuity of the operator $\sigma ^*\hat{P} \Upsilon
_{\sigma ^*}|_{S^{n-1}}$, Identities $(2)$, Theorem 8 and Theorem
6.12 \cite{axlerb}  we deduce that the limit
$$\lim_{n\to \infty } h_n = h= \sigma ^*\hat{P} (\Upsilon _{\sigma
^*}|_{S^{n-1}}y)$$ exists in $K_p(B,{\cal X})$ relative to the norm
17$(1)$, consequently, $K_p(B,{\cal X})$ is the left Banach module
over the quaternion skew field.
\par According to Proposition 1.c.8
and Theorem 1.c.9 in volume 2 of the book \cite{lindliorb} the
Banach space $L^p(S^{n-1},{\bf R})$ has an unconditional basis. Then
we choose an unconditional basis $\{ f_m: ~ m\in {\bf N} \} $ in
$L^p(S^{n-1},{\bf H})$. Therefore, $\{ f_me_k: ~ m\in {\bf N},
k=1,...,l \} $ is the unconditional basis in the left $\bf H$ module
$L^p(S^{n-1},{\cal X}).$ For each $f_m$ then we choose a particular
solution $g_m\in W^p_1(S^{n-1},{\cal X})$ of the equation $\sigma
^*g_m=f_m$ restricted to $S^{n-1}$. \par One can put $Q$ to be the
left Banach module over $\bf H$ which is the closure in the Sobolev
space $W^p_1(S^{n-1},{\cal X})$ of $l-span_{\bf H} \{ g_me_k: ~ m\in
{\bf N}, k=1,...,l \} $. Therefore, we deduce that $\sigma
^*Q=L^p(S^{n-1},{\cal X})$. From Theorem 6.12 \cite{axlerb} and
Formulas $(1,2)$ above it follows that $\{ \sigma ^*\hat{P}(g_me_k):
m, k \} $ is the unconditional basis in $K_p(B,{\cal X})$.

\par {\bf 19. Remark.} More concrete bases can be constructed with the help of Theorem 15
and Corollary 16 as it is outlined below. In the real Hilbert space
$L^2(S^{n-1},{\bf R})$ supplied with the standard scalar product
$$ \{ f, g \} := \int_{S^{n-1}} f(y)g(y)\psi (dy)$$ the subspace of harmonic
polynomials restricted on $S^{n-1}$ is dense, where
$f$ and $g$ are functions from $S^{n-1}$ into ${\bf R}$, while $f$
and $g\in L^2(S^{n-1},{\bf R})$. Moreover, the decomposition
$$L^2(S^{n-1},{\bf R})= \bigoplus_{m=0}^{\infty }{\cal
H}_m(S^{n-1})$$ is valid (see Theorem 5.8 \cite{axlerb}), where
${\cal H}_m({\bf R}^n)$ denotes the space of all harmonic real
homogeneous polynomials $P_m(x)$ of degree $m$ on ${\bf R}$, that is
$P_m(tx)=t^mP_m(x)$ for each $x\in {\bf R}^n$ and $t\in {\bf R}$,
$\Delta P_m(x)\equiv 0$, whilst the vector space ${\cal
H}_m(S^{n-1})$ is the restriction of ${\cal H}_m({\bf R}^n)$ to
$S^{n-1}$. According to Theorem 5.34 \cite{axlerb} the set $ \{
D^{\alpha }|x|^{2-n}: |\alpha |=m, \alpha _1= 0 \mbox{ or } 1 \}
|_{S^{n-1}}$ is a basis for the real vector space ${\cal
H}_m(S^{n-1})$. Therefore,  $l-span_{\bf H} \bigcup_{m=0}^{\infty }
\{ D^{\alpha }|x|^{2-n}: |\alpha |=m, \alpha _1= 0 \mbox{ or } 1 \}
|_{S^{n-1}}$ is dense in $C^1(S^{n-1},{\bf R})$, since
$C^1(S^{n-1},{\bf R})\subset L^2(S^{n-1},{\bf R})$ and ${\cal
H}_m(S^{n-1})\subset C^1(S^{n-1},{\bf R})$ for each $m$.

\par From the decomposition of the Banach space $C(S^{n-1},{\bf H})=
\bigoplus _{j=0}^3 C(S^{n-1},{\bf R})i_j$ it follows that the left $\bf H$ span of the set
$ \{ D^{\alpha }\sigma ^*_x|x|^{2-n}: |\alpha |=m, \alpha _1= 0 \mbox{ or } 1; m= 1, 2,... \} |_{S^{n-1}}$
is dense in $C(S^{n-1},{\bf H})$, since $D^{\alpha }$ and $\sigma ^*$ in Cartesian coordinates $x_1,...,x_n$
commute, where $x=(x_1,...,x_n)\in {\bf R}^n$, $D^{\alpha }= \partial ^{|\alpha |} /\partial x_1^{\alpha _1}
...\partial x_n^{\alpha _n}$, $\alpha =(\alpha _1,...,\alpha _n)$, $\alpha _k\in \{ 0, 1, 2,...\}$ for each
$k=1,...,n$, $|\alpha | := \alpha _1+...+\alpha _n$. Using the selection procedure of \S 15
and Corollary 16 one gets from this system a Schauder basis in $K(B,{\cal Y})$.
\par It is possible also to take as the starting point Schauder bases in $C^1(S^{n-1},{\bf R})$.
In several works (see \cite{cisdom72,ryllsm78,semadb} and references therein) Schauder bases were constructed in the Banach spaces $C^m([0,1]^k,{\bf R})$.
\par Using two chart atlas of the $C^{\infty }$ Riemann manifold $S^{n-1}$ and a Schauder basis in $C^m(D,{\bf R})$ one can construct
a Schauder basis in $C^m(S^{n-1},{\bf R}),$ where $D := \{ x: x\in
{\bf R}^{n-1}; |x|\le 1 \} $, $m\ge 1$. In $C(S^{n-1},{\bf R})$ a
Schauder basis exists due to Weierstrass theorem and in
$C(S^{n-1},{\bf H})$ according to Theorem 2.4.
\par Above Clifford algebras and modules were considered over the quaternion skew
field ${\bf H}$. One can also consider Clifford algebras $\cal{X}$
and modules $\cal{Y}$ both over ${\bf R}$ or ${\bf C}$ and construct
Schauder bases in $K(B,\cal{Y})$ and in $K_p(B,\cal{Y})$ with
$1<p<\infty $ over the field either ${\bf R}$ or ${\bf C}$
correspondingly using results of this paper.
\par Apart from the left ${\bf H}$-module $K_p(B,\cal{Y})$ with
$1<p<\infty $, the left ${\bf H}$-module $K(B,\cal{Y})$ is not
expected to have an unconditional basis, since $C(0,1)$ does not
have an unconditional basis and does not even embed in a space with
an unconditional basis (see page 2 in volume 2 of the book
\cite{lindliorb}).
\par The results of this paper can be used for a subsequent investigation
of kernels $ker (\sigma )$ of the Dirac operator $\sigma $,
integration and solution of partial differential equations.

\end{document}